\documentclass[11pt,oneside]{article}
\usepackage{amssymb,amsmath}
\usepackage{url}
\usepackage[english]{babel}
\usepackage{pstricks,pst-node}
\renewcommand{\le}{\leqslant}
\renewcommand{\ge}{\geqslant}
\newcommand{\bad}{\mathbf{Bad}}
\newcommand{\mad}{\mathbf{Mad}}
\newcommand{\badl}{\mathbf{Mad}}

\newcommand{\RR}{\mathbb{R}}
\newcommand{\ZZ}{\mathbb{Z}}

\newcommand{\QQ}{\mathbb{Q}}
\newcommand{\KK}{\mathbf{K}}

\newcommand{\NN}{\mathbb{N}}

\newcommand{\JJ}{\mathcal{J}}
\newcommand{\II}{\mathcal{I}}
\newcommand{\LLL}{\mathcal{L}}
\newcommand{\DDD}{\mathcal{D}}

\newcommand{\KKK}{\mathbf{K}}
\newcommand{\LKK}{\mathbf{LK}}

\newcommand{\TTT}{\mathcal{T}}
\newcommand{\RRR}{\mathbf{R}}
\newcommand{\R}{\mathcal{R}}

\newcommand{\vr}{\mathbf{r}}
\newcommand{\vs}{\mathbf{s}}

\newtheorem{lemma}{Lemma}
\newtheorem{theorem}{Theorem}
\newtheorem{proposition}{Proposition}

\newtheorem{corollary}{Corollary}

\newtheorem{MLconj}{Mixed Littlewood Conjecture (MLC)}

\newtheorem{Lconj}{Littlewood Conjecture (LC)}

\newtheorem{MDP}{Mass Distribution Principle}

\newcommand{\cC}{{\cal C}}

\newcommand{\cJ}{{\cal J}}

\begin{document}


\title{Multiplicatively  badly approximable numbers and \\  generalised Cantor sets}

\author{Dzmitry Badziahin \footnote{Research supported by EPSRC grants EP/E061613/1}
\\ {\small\sc York}
  \and  Sanju Velani\footnote{Research supported by EPSRC grants EP/E061613/1 and EP/F027028/1 }
\\ {\small\sc York}}

\date{\it Dedicated to Andrew Pollington  \\ on
hitting  57*}


\maketitle

\begin{abstract}
Let $p$ be a prime number. The  $p$-adic case of the Mixed
Littlewood Conjecture states that  $\liminf_{q\to \infty}
q\cdot|q|_p\cdot||q\alpha||=0$  for all $\alpha \in \RR$. We show
that with the additional factor of $\log q \log\log q$ the statement is false. Indeed,  our main result implies that the set of $\alpha$ for
which $ \liminf_{q\to\infty} q\cdot\log q\cdot\log\log q\cdot
|q|_p\cdot||q\alpha||>0$ is of full dimension.  The result is obtained as an application of a general framework for  Cantor  sets developed in this paper.

\end{abstract}

\section{Introduction}

The goal of this paper is simple enough.  It is an attempt to
address the question:

\begin{itemize}
\item[] {\em What are the  analogues of the classical set of badly approximable numbers within the multiplicative frameworks of Littlewood's conjecture and its mixed counterpart?}
\end{itemize}

\subsection{The classical setup and  the set $\bad$}

A classical result of Dirichlet states that for any real number
$\alpha$ there exist infinitely many $q\in\NN$ such that
$$
q||q\alpha||<1   \   .
$$
Here and throughout  $||\, . \, ||$ denotes the distance to the
nearest integer. In general the right hand side of the above
inequality cannot be replaced by an arbitrarily small constant.
Indeed a result of Jarn\'{\i}k \cite{Ja} and Besicovitch \cite{Best}
states that the set
$$
\bad:=\{\alpha\in \RR\;:\; \liminf_{q\to \infty} q||q \alpha ||>0\}
$$
of badly approximable numbers is of maximal Hausdorff dimension;  i.e.
$$
\dim \bad=1.
$$
For details regarding  Hausdorff  dimension the reader is referred to  \cite{falc}.
However, from a measure theoretic point of view  the classical
theorem of Khintchine \cite{Kh24} enables us to improve on the
global statement (a statement true for all numbers) of Dirichlet by
a logarithm. In particular, for $\lambda\ge 0$  let
$$
\bad^{\lambda}  :=   \{\alpha\in \RR\;:\; \liminf_{q\to\infty} q
\cdot (\log q)^\lambda  \cdot ||q\alpha||>0\} \ .
$$
Then, Khintchine's theorem implies that
$$
\big|\bad^{\lambda} \big|=\left\{\begin{array}{ll}
0& \mbox{ if }\lambda\le 1\\[2ex]
{\mbox{\footnotesize FULL}}& \mbox{ if }\lambda> 1.
\end{array}\right.
$$
Here and throughout  $|\, \cdot \, |$ denotes Lebesgue measure and
`{\footnotesize FULL}' means that the complement of the set under
consideration is of measure zero.

The upshot of the classical setup is that we are able to shave off a
logarithm from the measure theoretic `switch over' set $\bad^1$
before  we precisely hit  the set  $\bad$. In addition, if we shave
off any more (i.e. $(\log q)^{1 + \epsilon}$ with $\epsilon > 0$
arbitrary) then the  corresponding set  becomes empty. This is  a
theme which we claim reoccurs  within  the  multiplicative framework
of Littlewood's conjecture and its mixed counterpart.

%
%

\subsection{The multiplicative setup and  the set $\mad$}

A straightforward consequence of Dirichlet's classical result is
that for every  $(\alpha,\beta) \in \RR^2$, there exist infinitely
many $q\in\NN$ such that
$$
q\cdot ||q\alpha||\cdot||q\beta||<1.
$$
Littlewood conjectured  that the right hand side  of the above
inequality can be replaced by an arbitrarily small constant.
\begin{Lconj}
 For every  $(\alpha,\beta) \in \RR^2$,
\begin{equation}\label{eq_LC}
\liminf_{q\to \infty}q\cdot ||q\alpha||\cdot ||q\beta||=0  \ .
\end{equation}
\end{Lconj}

\noindent Despite concerted efforts over the years this famous
conjecture remains open.   For background and recent `progress'
concerning  this fundamental  problem  see \cite{ELK,PVL} and
references therein.

A consequence of LC   is that the set
$$
\{(\alpha,\beta)\in \RR^2\;:\; \liminf_{q\to\infty} q\cdot
||q\alpha||\cdot ||q\beta||>0\}
$$
is empty and therefore is not a candidate for the multiplicative
analogue of $\bad$. Regarding possible candidates, for $\lambda\ge
0$ let
$$
\mad^{\lambda}:=\{(\alpha,\beta)\in \RR^2\;:\; \liminf_{q\to\infty}
q \cdot (\log q )^\lambda \cdot ||q\alpha||\cdot ||q\beta||>0\}.
$$
From a measure theoretic point of view Gallagher's theorem \cite{Ga}
(the multiplicative analogue of Khintchine's theorem) implies that
$$
|\mad^\lambda|=\left\{\begin{array}{ll}
0& \mbox{ if }\lambda\le 2\\[2ex]
{\mbox{\footnotesize FULL}}& \mbox{ if }\lambda> 2.
\end{array}\right.
$$
Natural  heuristic `volume' arguments
give evidence in favour of the following statement: for every
$(\alpha,\beta) \in \RR^2$ there exist infinitely many $q\in\NN$
such that
$$
q\cdot \log q \cdot  ||q\alpha||\cdot||q\beta|| \ll 1  .
$$
The  results of  Peck \cite{Peck} and Pollington $\&$ Velani
\cite{PVL}  give solid support to this statement which represents a
significant strengthening of  Littlewood's conjecture and    implies
that

\begin{itemize}
\item[]  {\bf[L1]} \hspace*{1ex}   $\badl^{\lambda}    =  \emptyset $ \ if \  $ \lambda < 1 $.
\end{itemize}
Moreover, we suspect that the heuristics are sharp  and thus  $\mad
:= \mad^{1} $ represents the natural analogue of $\bad$ within the
multiplicative setup.  It is worth emphasizing that $\mad$ defined
in this manner is precisely the set we hit after shaving  off a
logarithm from the measure theoretic `switch over' set $\mad^2$.
Note that  this is in keeping with the classical setup. Furthermore,
we claim that the analogue of Jarn\'{\i}k-Besicovitch  theorem is
true for $\mad$. In other words,
\begin{itemize}
\item[]  {\bf[L2]} \hspace*{1ex}  $ \dim \badl^{\lambda}    =  2 $ \ if  \  $ \lambda  \ge  1 $.
\end{itemize}

\vspace*{2ex}

\noindent

 Regarding  {\bf [L1]}, notice that  a counterexample to  LC would imply that $\badl^{\lambda} $ is non-empty    for any  $\lambda \ge 0$.  In principle,  it should be easier to give a counterexample to {\bf[L1]}. To date all that is know is the remarkable  result of  Einsiedler, Katok $\&$ Lindenstrauss \cite{ELK} that states that
$ \dim \mad^0=0. $ The following would be a leap in the right
direction towards {\bf[L1]} and  would  represent a significant
strengthening of the Einsiedler-Katok-Lindenstrauss zero dimension
result.

\begin{itemize}
\item[]  {\bf[L3]} \hspace*{1ex}
$\dim \mad^\lambda =0  \ \mbox{ if }  \  \lambda<1.$
\end{itemize}

\noindent To the best of our knowledge, currently we do not even
know if  $ \dim \badl^{\lambda}    < 2 $  for  strictly positive $
\lambda < 1 $.

\noindent

 Regarding {\bf [L2]}, very little beyond the trivial is known. A simple consequence of the `{\mbox{\footnotesize FULL}}'  statement above is that $ \dim  \mad^\lambda =  2  $ if $\lambda > 2$.  Recently, Bugeaud $\&$ Moshchevitin~\cite{BM} have shown that  $\dim \mad^2=2$. Note that this is non-trivial since the set $ \mad^2$  is of measure zero.  Surprisingly and somewhat embarrassingly we are unable to show that  $ \mad^{2 - \epsilon} \neq \emptyset $   let alone

\begin{itemize}
\item[]  {\bf[L4]} \hspace*{1ex}
$ \mad^\lambda    \neq \emptyset   \ \mbox{ if }  \  1 \le
\lambda<2.$
\end{itemize}

\noindent In other words, given our current state of knowledge, we
can not rule out the unlikely  possibility that LC is actually true
with a  $(\log q)^{2- \epsilon}$ term  inserted in the left hand
side of~\eqref{eq_LC} -- see also \cite[Question 37]{aim}.

\vspace*{2ex}

In this paper, we are unable to  directly contribute  towards the
statements  {\bf[L1]} -- {\bf[L4]}.   However,  we are able to make
a significant contribution towards  establishing the analogue of
{\bf[L2]} within the framework of the mixed Littlewood  conjecture.
Thus,  if there is a genuine `dictionary' between the results
related to the two conjectures  then indirectly  our contribution
adds  weight  towards {\bf[L2]}.

\subsection{The mixed multiplicative setup and the set $\badl_{\DDD}$}

Recently, de Mathan   $\&$  Teuli\'{e} in  \cite{MT04} proposed the
following variant of Littlewood's conjecture.  Let $\mathcal{D} $ be
a sequence $(d_k)_{k=1}^\infty$ of integers greater than or equal to
$2$ and  let
$$
 D_0 := 1     \quad {\rm and}  \quad    D_n := \prod_{k=1}^n d_k   \ .
$$
For  $q \in \ZZ$ set
$$
|q|_\DDD:=\inf\{D_n^{-1}\;:\; q\in D_n\ZZ\}.
$$

\begin{MLconj}
 For every real number $\alpha$
\begin{equation}\label{eq_MLC}
\liminf_{q\to \infty}q\cdot |q|_\DDD\cdot ||q\alpha||=0  \ .
\end{equation}
\end{MLconj}

\noindent When $\mathcal{D}$ is the constant sequence equal to a
prime number $p$,  the norm $|\, \cdot \, |_\mathcal{D}$ is the
usual $p$-adic norm $|\, \cdot \, |_p$.   In this particular case,
there is a perfect dictionary   between the current body of results
associated with ($p$-adic) MLC and LC. The following constitute the
main non-trivial entries.

\begin{itemize}
\item In  \cite[Theorem 2.1]{MT04}  de Mathan   $\&$  Teuli\'{e} establish  the analogue of Peck's cubic result.
\item In \cite[Section 1]{MT04}  de Mathan   $\&$  Teuli\'{e}  observe that the ideas within \cite{PVL} establish the analogue of the Pollington-Velani full dimension result. Also see \cite[Theorem 4]{bdd}.
\item In \cite[Theorem 1]{BHV} Bugeaud, Haynes $\&$ Velani establish the analogue of  Gallagher's measure theoretic  result.
\item In \cite[Theorem 1.1]{EinsiedlerKleinbock} Einsiedler $\&$ Kleinbock establish the analogue of the Einsiedler-Katok-Lindenstrauss zero dimension  result.
\item In \cite{BM} Bugeaud $\&$ Moshchevitin establish the analogue of their  $\dim \mad^2=2$ result.
\end{itemize}

\noindent Moving away from the $p$-adic case, the results associated
with MLC in the first two items above are valid for any bounded
sequence $\DDD$. In all likelihood, this is also true for the other
three items. The biggest  challenge  of the three seems to lie in
generalising  the ($p$-adic)  result of  Einsiedler $\&$ Kleinbock
to bounded sequences.  We are pretty confident that the other two
items can be generalised to bounded $\DDD$  without too much trouble
but stress that we have not carried out the details\footnote{The
problem of generalizing the ($p$-adic) mixed result obtained in
\cite{BHV}  to arbitrary sequences $\DDD$ is particularly
interesting  since for unbounded $\DDD$ we suspect that the `volume'
sum is dependant on $\DDD$.}. The point being made here is that for
bounded $\DDD$ there is reasonably hard  evidence in support of  a
`LC--MLC' dictionary.

For $\lambda \ge 0$ let
\begin{equation}
\badl_{\DDD}^{\lambda} := \left\{\alpha\in \RR \;:\; \liminf_{q \to
\infty} \  q \cdot (\log q)^{\lambda}
 \cdot
|q|_\DDD\cdot ||q\alpha||> 0   \right\}.
\end{equation}

\noindent For  bounded $\DDD$, in view of the  above discussion it
is natural to expect that the following statements correspond to the
entries  {\bf[L1]} and {\bf[L2]} within the `LC--MLC' dictionary.


\begin{itemize}
\item[]  {\bf[ML1]} \hspace*{1ex}   $\badl_{\DDD}^{\lambda}    =  \emptyset $ if $ \lambda < 1 $.

\item[]  {\bf[ML2]} \hspace*{1ex}  $ \dim \badl_{\DDD}^{\lambda}    =  1 $ if $ \lambda  \ge  1 $.
\end{itemize}

\noindent In short,  the upshot for bounded $\DDD$ is that
$\badl_{\DDD} := \badl_{\DDD}^{1} $ represents the natural analogue
of $\bad$ within the `mixed' multiplicative setup. The assumption
that $\DDD$ is bounded is absolutely necessary -- see Theorem
\ref{th_mainsvsv}   below.

Obviously a counterexample to MLC would imply that
$\badl_{\DDD}^{\lambda}    \neq   \emptyset$   for any  $\lambda
\ge 0$.  In principle,  it should be easier to give a
counterexample to {\bf[ML1]}. The  Einsiedler--Kleinbock  result
($\dim \badl_{\DDD}^{0} = 0$  within the $p$-adic case) represents
the current state of knowledge regarding {\bf[ML1]}. It would be
highly desirable to obtaining the following generalization.
\begin{itemize}
\item[]  {\bf[ML3]} \hspace*{1ex}   $\dim \badl_{\DDD}^{\lambda}    =  0  $ if $ \lambda < 1 $.
\end{itemize}
\noindent As far as we are aware,  it is not  even known  if  $ \dim
\badl_{\DDD}^{\lambda}    < 1 $  for  strictly positive $ \lambda <
1 $.

The following contribution towards  {\bf[ML2]}  constitutes the main
result proved in this paper.  In our opinion, up to  powers of
logarithms it is best possible for bounded $\DDD$.

\begin{theorem}\label{th_main}
Let $\DDD$ be a  sequence of integers greater than or equal to $2$.
Then the set of real numbers $\alpha$ such that
\begin{equation}\label{eq_main}
\liminf_{q\to\infty}\   q\cdot \log q\cdot \log\log q\cdot
|q|_\DDD\cdot ||q\alpha||>0.
\end{equation}
has Hausdorff dimension equal to 1.
\end{theorem}

A simple consequence of the theorem is  the following statement.

\begin{corollary}\label{mainsv1}
Let $\DDD$ be a  sequence of integers greater than or equal to $2$.
For  $ \lambda  >   1  $
$$
\dim \badl_{\DDD}^{\lambda}    =  1     \, .
$$
\end{corollary}

\noindent Unfortunately, for bounded $\DDD$ we are unable to deal
with the case $\lambda =1$. In fact, we are unable to show that

\begin{itemize}
\item[]  {\bf[ML4]} \hspace*{1ex}
$\badl_{\DDD}^{1} \neq \emptyset  \, . $
\end{itemize}

\noindent However, for unbounded $\DDD$  we can do much better in
the following sense.

\begin{theorem}\label{th_main2}
Let $\DDD=\{2^{2^n}\}_{n\in\NN}$. Then the set of real numbers
$\alpha$ such that
\begin{equation}\label{eq_dunbound}
\liminf_{n\to\infty} q\cdot \log\log q\cdot \log\log\log q\cdot
|q|_\DDD\cdot ||q\alpha||>0
\end{equation}
has Hausdorff dimension equal to $1$.
\end{theorem}

A simple consequence of the theorem is  the following statement.

\begin{corollary}\label{th_mainsvsv}  There exist uncountably many unbounded  sequences $\DDD$ of
integers greater than or equal to $2$ such that
\begin{equation}\label{eq_mainsvsv}
\dim \badl_{\DDD}^{\lambda}    =  1   \quad \forall \  \lambda > 0
\ .
\end{equation}
\end{corollary}

\noindent  The theorem shows  that  {\bf[ML1]} is not generally true
for unbounded $\DDD$. It also suggests that if there are
counterexamples to MLC then they may be easier to find among rapidly
increasing  sequences.
 Furthermore, for unbounded $\DDD$ it is not generally true that the  natural analogue of $\bad$ within the `mixed' multiplicative setup is  $\badl_{\DDD}^{1} $. This is yet an other reason to why we restrict the `LC--MLC' dictionary to bounded sequences.  Indeed, we can deduce from the proof of Theorem~\ref{th_main2}  that the analogue of $\bad$ for any given  unbounded $\DDD$  is in fact  dependant on the growth of~$\DDD$.


\section{Preliminaries}


To prove Theorems~\ref{th_main} and~\ref{th_main2} it will be
convenient to work with the `modified logarithm' function
$\log^*\;:\; \RR\to\RR$ defined as follows
$$
\log^*\! x:=\left\{\begin{array}{ll}1&\mbox{for
}x<e \\
\log x&\mbox{for }x\ge e \, .
\end{array}\right.
$$
This will guarantee that for small values of $x$ the function
$\log^*\!x$ is well defined.

\subsection{The basic strategy \label{bssv}}

Given a function $f\;:\;\NN\to \RR$ and a sequence $\DDD$  of
integers not smaller than 2, consider the set
\begin{equation}\label{def_badl}
\mad_\DDD(f):=\{\alpha\in \RR\;:\; \liminf_{q\to \infty} f(q)\cdot
q\cdot |q|_\DDD\cdot ||q\alpha||>0\}.
\end{equation}
By definition the set $\mad_\DDD(f)$ is a subset of $\RR$  and
therefore
 $$\dim \mad_\DDD(f) \le 1.$$ Thus the proofs of Theorems~\ref{th_main} and~\ref{th_main2} are reduced to establishing the following respective  statements.
\begin{proposition}\label{prop1sv}
Let $\DDD$ be a  sequence of integers greater than or equal to $2$.
Then
\begin{equation}\label{eq_madf}
\dim\mad_\DDD(f)\ge 1\quad\mbox{with }  \ f(q):=\log^*\!q\cdot\log^*\log q.
\end{equation}
\end{proposition}
\begin{proposition}\label{prop2sv}
Let $\DDD=\{2^{2^n}\}_{n\in\NN}$. Then
\begin{equation}\label{eq_madf2}
\dim\mad_\DDD(f)\ge 1\quad\mbox{with }  \ f(q):=\log^*\log
q\cdot\log^*\log^*\log q.
\end{equation}
\end{proposition}

To establish  the propositions  we make use of the following
decomposition of $\mad_\DDD(f)$. For any constant $c>0$ define
$$
\badl_{\DDD}(f,c) := \left\{\alpha\in \RR \;:\; f(q)\cdot q\cdot
|q|_\DDD\cdot ||q\alpha||>c\;\ \forall   \  q\in\NN\right\}.
$$
It is easily verified that $$ \badl_{\DDD}(f,c) \subset
\badl_\DDD(f)$$  and
$$\badl_\DDD(f)\, =   \, \bigcup_{c > 0}  \badl_{\DDD}(f,c)\ . $$

\noindent Geometrically, the set $\badl_{\DDD}(f,c) $ simply
consists of points on the real line that avoid all intervals
$$
\Delta(r/q):=\left[\frac{r}{q}-\frac{c}{f(q)q^2|q|_\DDD},\frac{r}{q}+\frac{c}{f(q)q^2|q|_\DDD}\right]
$$
centered at rational points $r/q $   with $q \ge 1$. Alternatively, points
on the real line that lie within any such interval are removed.
Given a rational  $r/q$,  let
\begin{equation}\label{def_height}
H(q):=q^2|q|_\DDD
\end{equation}
denote its {\em height}.  Trivially,   we have that
$$
|\Delta(r/q)|=\frac{2c}{f(q)H(q)}.
$$

\noindent  In order to show that   $\dim \badl_{\DDD}(f) \ge 1  $,
the idea is to construct a Cantor-type subset $\KK_{\DDD}(f,c)   $
of $
 \badl_{\DDD}(f,c) $   such that
$$
\dim \KK_{\DDD}(f,c)    \ge 1
$$
for some small constant $c>0$. Hence, by construction we have that
$$
\dim\mad_\DDD(f)   \, \ge \, \dim \badl_{\DDD}(f,c) \, \ge \, \dim
\KK_{\DDD}(f,c)   \, \ge \,   1\, .
$$
Thus, the name of the game is to construct the `right type' of Cantor
set $\KK_{\DDD}(f,c)$.   In short, the  properties of the desired
set  fall naturally within a general framework which we now
describe.

\subsection{A general Cantor  framework \label{gcf}}

\noindent{\bf The parameters.}  Let ${\rm I}$ be a closed interval
in $\RR$. Let $$\RRR:=(R_n)  \quad {\rm with }  \quad {n\in \ZZ_{\ge
0}}$$ be a sequence of natural numbers and $$\vr:=(r_{m,n})   \quad
{\rm with }  \quad  m,n\in \ZZ_{\ge 0} \ {\rm  \ and \  }  \ m\le n
$$ be a two parameter sequence of non-negative real numbers.

\vspace*{2ex}

\noindent{\bf The  construction.}  We start by subdividing the interval ${\rm I}$ into $R_0$
closed intervals $I_1$  of equal length and denote by $\II_1$ the collection of such intervals.  Thus,
$$
\#\II_1  =  R_0   \qquad {\rm and } \qquad |I_1| =  R_0^{-1}\, |{\rm
I}|  \ .
$$
Next,  we remove  at most  $r_{0,0}$ intervals $I_1$ from $\II_1$ .   Note that we do not specify which intervals should be
removed but just give an upper bound on the number of intervals to be removed. Denote by  $\JJ_1$ the resulting
collection. Thus,
\begin{equation}\label{iona1}
\#\JJ_1  \ge    \#\II_1   -  r_{0,0}  \, .
\end{equation}
For obvious reasons, intervals in $\JJ_1$ will be referred to as (level one) survivors.   It will be convenient to define  $\JJ_0 := \{J_0\} $ with $ J_0 :={\rm I} $.

\vspace*{1ex}

\noindent In general, for $n \ge 0$, given  a collection
$\JJ_n$    we construct a nested collection $\JJ_{n+1}$  of closed intervals $J_{n+1}$  using the following two operations.
\begin{itemize}
\item{\em Splitting procedure.} We subdivide each  interval $J_n\in \JJ_n$  into $R_n$ closed sub-intervals $I_{n+1}$ of equal length and denote by  $\II_{n+1}$ the collection of such intervals. Thus,
    $$
    \#\II_{n+1}  =  R_n \times  \#\JJ_n   \qquad {\rm and } \qquad |I_{n+1}| = R_n^{-1}\, |J_n|    \ .
    $$
\item{\em Removing procedure.} For each interval $J_n\in \JJ_n$ we remove at most
$r_{n,n}$ intervals $I_{n+1} \in  \II_{n+1} $ that lie  within $
J_n$.  Note that  the number of intervals $I_{n+1}$ removed is allowed to  vary amongst  the  intervals  in  $\JJ_n$.    Let  $\II_{n+1}^{n} \subseteq \II_{n+1}  $ be the collection of  intervals that remain.  Next, for each interval $J_{n-1}\in \JJ_{n-1}$
we remove at most  $r_{n-1,n}$  intervals $I_{n+1} \in  \II_{n+1}^{n} $ that lie within $ J_{n-1}$. Let  $\II_{n+1}^{n-1} \subseteq \II_{n+1}^{n}  $ be the collection of  intervals that remain. In general, for each interval $J_{n-k}\in \JJ_{n-k}$  $(1 \le k \le n)$
we remove  at most $r_{n-k,n}$ intervals $I_{n+1} \in \II_{n+1}^{n-k+1} $ that lie within $J_{n-k}$. Also we let  $\II_{n+1}^{n-k} \subseteq \II_{n+1}^{n-k+1}  $ be the collection of  intervals that remain.
In particular,
$\JJ_{n+1}  := \II_{n+1}^{0} $
is the desired collection of (level $n+1$) survivors.  Thus,
the total number of intervals $I_{n+1}$ removed during the removal procedure  is at most
$
r_{n,n}\#\JJ_n+r_{n-1,n}\#\JJ_{n-1}+\ldots+r_{0,n}\#\JJ_0
$
and so
\begin{equation}\label{iona2}
\#\JJ_{n+1}\ge R_n\#\JJ_n-\sum_{k=0}^nr_{k,n}\#\JJ_k.
\end{equation}
\end{itemize}
\noindent Finally, having constructed the nested collections $\JJ_n$ of closed intervals   we consider the limit set
$$
 \KKK ({\rm I},\RRR,\vr) :=  \bigcap_{n=1}^\infty \bigcup_{J\in
\JJ_n} J.
$$
The set $\KKK({\rm I},\RRR,\vr)$ will be referred to as a {\em $({\rm I},\RRR,\vr)$ Cantor set.}

\vspace*{2ex}
\noindent{\em Remark.}
We stress  that the  triple $({\rm I},\RRR,\vr)$ does not
uniquely determine the set $\KKK ({\rm I},\RRR,\vr)$.
The point is that during the construction  we only specify the maximum number of intervals rather than the specific intervals  to be removed. Thus the triple $({\rm I},\RRR,\vr)$ gives rise to a family of $({\rm I},\RRR,\vr)$ Cantor sets that reflects the various available choices  during the removing procedure.

\vspace*{1ex}
As an illustration of the general framework, it is easily seen that the standard middle third
Cantor set   corresponds to a $({\rm I}, \RRR,\vr)$ Cantor set  with
$$
{\rm I} := [0,1], \qquad \RRR= (3,3,3,\ldots)\quad\mbox{and}\quad
\vr=(r_{m,n})
$$
where
$$
r_{m,n}:=\left\{\begin{array}{ll}1&\mbox{if }m=n\\[1ex]
0&\mbox{otherwise}.
\end{array}\right.
$$
\vspace*{2ex}

\noindent{\bf The results.} By definition, if $\JJ_n$ is empty for some $n \in \NN$ then the corresponding set  $ \KKK ({\rm I},\RRR,\vr)$ is obviously empty.    On the other hand, by construction,  each closed interval $J_n \in \JJ_n $ is contained in some closed interval   $J_{n-1} \in \JJ_{n-1} $. Therefore
 $$\KKK({\rm I}, \RRR,\vr) \neq \emptyset   \qquad {\rm if }  \qquad  \# \JJ_n \ge 1  \quad \forall  \  n \in \NN \, .$$
Our first result provides a natural condition that  guarantees this cardinality hypothesis   and  therefore the non-empty statement.

\begin{theorem}\label{th_cantor1} Given $\KKK ({\rm I},\RRR,\vr) $,
let
\begin{equation}\label{def_t0}
t_0:=R_0-r_{0,0}
\end{equation}
and for $n \ge 1 $ let
\begin{equation}\label{def_tn}
t_n:=R_n-r_{n,n}-\sum_{k=1}^n \frac{r_{n-k,n}}{\prod_{i=1}^k
t_{n-i}}    \,    .
\end{equation}
Suppose that  $t_n>0$ for all $n\in\ZZ_{\ge 0}$. Then
$$
\KKK ({\rm I},\RRR,\vr)  \neq \emptyset \ .
$$
\end{theorem}

\vspace*{2ex}

The proof of  Theorem~\ref{th_cantor1} is short and direct and there seems little point in delaying it.

\vspace*{1ex}

\noindent{\em Proof of Theorem~\ref{th_cantor1}.} We show that a  consequence of the construction of $\KKK ({\rm I},\RRR,\vr) $ is that
\begin{equation}\label{ineq_jj}
\#\JJ_n\ge t_{n-1} \#\JJ_{n-1} \quad \forall  \  n \in \NN \,  .
\end{equation}
This together with the assumption that  $t_n>0$  implies that
$
\#\JJ_n\ge \prod_{i=0}^{n-1} t_i \, \#\JJ_0  \, > \, 0
$
and thereby completes  the proof of Theorem~\ref{th_cantor1}.  To verify \eqref{ineq_jj} we use induction.  In view of \eqref{iona1} and \eqref{def_t0} the statement is trivially true for $n=1$.
Now suppose that \eqref{ineq_jj} is true for all $1 \le k \le n$. In particular, for any  such $k$ we have that
$$
\#\JJ_n\ge t_{n-1}\#\JJ_{n-1}\ge\ldots\ge \prod_{i=1}^{k}t_{n-i}\#
\JJ_{n-k}.
$$
\noindent Thus,
\begin{eqnarray*}
\#\JJ_{n+1}  & \stackrel{\eqref{iona2}} \ge  &  R_n\#\JJ_n-\sum_{k=0}^nr_{n-k,n}\#\JJ_{n-k}
\\ &\ge &
R_n\#\JJ_n-r_{n,n}\#\JJ_n-\sum_{k=1}^n\frac{r_{n-k,n}\#\JJ_n}{\prod_{i=1}^k
t_{n-i}}  \\
& \stackrel{\eqref{def_tn}}  =  &  t_n\#\JJ_n.
\end{eqnarray*}
This completes the induction step and establishes  \eqref{ineq_jj} as required.
\hspace*{\fill}$\boxtimes$

\vspace*{3ex}

Our next result  enables us to estimate the Hausdorff dimension of
$\KKK({\rm I},\RRR,\vr)$.  It is the key to establishing Propositions \ref{prop1sv} $\&$ \ref{prop2sv}.

\begin{theorem}\label{th_cantor2}  Given $\KKK ({\rm I},\RRR,\vr) $, suppose that
 $R_n\ge 4$ for all $n\in\ZZ_{\ge 0}$ and  that
\begin{equation}\label{cond_th2}
\sum_{k=0}^n \left(r_{n-k,n}\prod_{i=1}^k
\left(\frac{4}{R_{n-i}}\right)\right)\le \frac{R_n}{4}.
\end{equation}
Then
$$
\dim \KKK ({\rm I},\RRR,\vr)  \ge \liminf_{n\to \infty}(1-\log_{R_n} 2).
$$
\end{theorem}

\noindent Here we use the convention that the product term in \eqref{cond_th2} is one when $k=0$  and by definition  $ \log_{R_n}\!2 := \log2/ \log R_n$.  The proof of Theorem~\ref{th_cantor2} is rather involved and constitutes the main substance of \S\ref{Stalin} and \S\ref{Lenin}. To some extent the raw ideas required to establish Theorem~\ref{th_cantor2} can be found in \cite[\S7]{BPV} where a conjecture of  W.M. Schmidt regarding the intersection of simultaneously badly approximable sets is proved.  Nevertheless we stress that in this paper we develop a general Cantor type framework rather than address a specific problem.
As a consequence the key ideas of  \cite{BPV} are foregrounded.


\vspace*{2ex}

\noindent {\em Remark.} Although Theorem~\ref{th_cantor2} is more than sufficient for  the specific application we have in mind,  we would like to point out  that we have not attempted to establish the most general or best possible statement.  For example, in the case $ R_n \to \infty$ as $n \to \infty$, the theorem  together with the fact that $ \KKK ({\rm I},\RRR,\vr) \subset \RR$ implies that  $\dim \KKK ({\rm I},\RRR,\vr) = 1   $.
However, we do not claim that condition  \eqref{cond_th2}  is optimal for establishing this full dimension result.

\vspace*{2ex}

In the final section of the paper, we show that the intersection of any finite number of sets $\KKK({\rm I},\RRR,\vr_i) $ is yet another  $({\rm I},\RRR,\vr) $ Cantor set for some  appropriately chosen  $\vr$. We shall also see in \S\ref{FC} that this enables  us to strengthen Theorem~\ref{th_main}.

\section{Proof of Proposition \ref{prop1sv} modulo Theorem~\ref{th_cantor2}   \label{proofcone}}

Throughout,  $f:\NN\to \RR  :  q \to f(q)=\log^*\!q\cdot\log^*\!\log q$
and $\DDD$  is  a  sequence of integers greater than or equal to
$2$.
Let
$$R > e^{12} $$
be an integer. Choose $c_1=c_1(R) > 0 $
sufficiently small so that
\begin{equation}\label{ineq_c1}
2e^2c_1\frac{\log R+2}{\log 2}R<1
\end{equation}
and let $c = c(R,c_1)>0$ be a constant such that
\begin{equation}\label{ineq_c}
c\left(\frac{64R^2(\log R+2)}{c_1\log 2}+\frac{16eR^2(\log
R+2)^2}{\log 2}\right)<1.
\end{equation}

With reference to \S\ref{bssv}  we  now describe the basic   construction of the  set $\KK_{\DDD}(f,c)$.
Let ${\rm I}$ be any interval of length $c_1$ contained within the unit interval $ [0,1]$.   Denote by  $ \JJ_0:=\{J_0\}$ where $J_0 := {\rm I}$.
The idea is to establish,  by induction on  $n$, the
existence of a collection $\JJ_n$ of closed intervals $J_n$ such
that $\JJ_n $ is nested in $\JJ_{n-1}$;   that is, each interval
$J_n$ in $\JJ_n$ is contained in some interval $J_{n-1}$ in
$\JJ_{n-1}$.    The length of an interval  $J_n$  will be given by
$$ |J_n|  \, :=  \,   c_1    \, R^{-n}F^{-1}(n)   \  ,
$$
where
$$
F(n):=\prod_{k=1}^n k \; [\log^*\! k]\;\  \mbox{ for } n\ge 1\quad\mbox{
and }  \quad F(0):=1  \   \mbox{ for }n\le 0  .
$$
Moreover, each  interval $J_n $ in $\JJ_n$ will satisfy the condition that
\begin{equation}  \label{cond}
J_{n} \,  \cap \, \Delta(r/q) \, =  \, \emptyset   \qquad    \forall
\ \  r/q\in\QQ  \ \  \mbox{with } \  H(q) <  R^{n \!-\!1}F(n-1) \, .
\end{equation}
In particular, we put
$$\KK_{\DDD}(f,c) :=   \bigcap_{n=1}^\infty \bigcup_{J\in\JJ_n}J \ .
$$
By construction, condition \eqref{cond} ensures that
$$
\KK_{\DDD}(f,c)  \subset \badl_{\DDD}(f,c) \ .
$$
Furthermore,  with reference to \S\ref{gcf} it will be apparent from the construction of the collections $\JJ_n $ that
$ \KK_{\DDD}(f,c) $  is in fact a  $({\rm I},\RRR,\vr)$ Cantor set  with   $\RRR=(R_n)$ given by
\begin{equation}\label{def_rn}
R_n:=R \, (n+1) \, [\log^*\!(n+1)]
\end{equation}
and  $\vr=(r_{m,n})$  given by
\begin{equation}\label{def_rrnm}
r_{m,n} \, :=  \, \left\{\begin{array}{ll} 7 \, \log^2\!R\cdot
n^2(\log^*\!n)^2 &\mbox{
\ if }\; m=n-1\\[2ex]
0    &\mbox{ \ otherwise. }
\end{array}\right.
\end{equation}

\noindent By definition, note that   for  any
$ R > e^9$ we have that the
\begin{eqnarray*}
{\rm l.h.s. \ of \ } \eqref{cond_th2}  & = &  r_{n-1,n}\cdot\frac{4}{R_{n-1}}  \, \le  \,    7 \cdot  2^3  \cdot \frac{\log^2\!R \cdot n \, \log^*\!n }{R}  \\[2ex]
&\le  &      \frac{  7 \cdot  2^6  \log^2\!R}{R^2} \cdot  \frac{R \, (n+1) \, [\log^*\!(n+1)] }{4} \\[2ex]
& \le &  \frac{R_n  }{4}  \ =  \  {\rm r.h.s. \ of \ }  \eqref{cond_th2}  \, .
\end{eqnarray*}
Since we are assuming that  $ R > e^{12}  $,  it then follows via   Theorem~\ref{th_cantor2}   that
$$
\dim \KK_{\DDD}(f,c)   \ge \liminf_{n\to \infty} (1-\log_{R_n}\!2)=1   \, .
$$
This  completes the proof of Proposition \ref{prop1sv}  modulo Theorem~\ref{th_cantor2} and the construction of  the collections $\JJ_n $.

%

\subsection{Constructing the collections $\JJ_n$ \label{cone}}


For $n=0$, we trivially have that (\ref{cond}) is satisfied for the
interval $J_0 = {\rm I}$. The point is that  there are no rationals satisfying
the height condition  $H(q)<1$ since by definition  $H(q) \ge 1$. For the same reason
(\ref{cond}) with $n=1$ is trivially satisfied for any interval
$J_1$  obtained by subdividing each $J_0$ in $\JJ_0$ into $R_0=R$
closed intervals of equal length $c_1 R^{-1} $. Denote by $\JJ_1 $
the resulting collection of intervals $J_1$ and note that
$
\# \JJ_1    =   R  \ .
$

In general, given $\JJ_n$ satisfying \eqref{cond} we wish to
construct a nested collection $\JJ_{n+1}$ of intervals $J_{n+1}$ for
which (\ref{cond}) is satisfied with $n$ replaced by $n+1$. By
definition, any interval $J_n$ in $\JJ_n$ avoids intervals
$\Delta(r/q)$ arising from rationals with height bounded above by the quantity
$R^{n-1}F(n-1)$. Since any `new' interval $J_{n+1}$ is to be nested
in some $J_n$, it is enough to show that $J_{n+1}$ avoids intervals
$\Delta(r/q)$ arising from rationals $r/q$ with height satisfying
\begin{equation}\label{zeq2}
R^{n-1}F(n-1)\le H(q)<R^nF(n)  \   .
\end{equation}
Denote by $C(n)$ the collection of all rationals satisfying
this height condition.  Formally
$$
C(n) := \left\{r/q \in \QQ  \, : \,  H(q)  \ \ {\rm satisfies \
(\ref{zeq2}) \, } \right\}   \
$$
and it is precisely this  collection of rationals that comes into
play when attempting to construct   $\JJ_{n+1}$ from $\JJ_{n}$. We now proceed
with the construction.

Assume that $n\ge 1$. We subdivide each $J_n$ in $\JJ_n$  into
$$
R_n\stackrel{\eqref{def_rn}}=R(n+1) \, [\log^*\!(n+1)]
$$
closed intervals $I_{n+1}$ of equal length $ c_1
R^{-(n+1)}F^{-1}(n+1)$ and denote by $\II_{n+1}$ the collection of
such intervals. Thus,
$$|I_{n+1}|=c_1 R^{-(n+1)}F^{-1}(n+1)    $$
and
$$
\# \II_{n+1}  \, =  \,  R(n+1) \ [\log^*\!(n+1)]  \, \times \, \# \JJ_{n}
\, .
$$
It is obvious that the construction of $\II_{n+1}$ corresponds
to the splitting procedure associated with the  construction of a  $({\rm I},\RRR,\vr)$ Cantor set.

In view of the nested requirement, the collection $\JJ_{n+1}$ which
we are attempting to construct will be a  sub-collection  of
$\II_{n+1}$. In other words, the intervals $I_{n+1}$ represent
possible candidates for $J_{n+1}$. The goal now is simple --- it is
to remove those `bad' intervals $I_{n+1}$ from $\II_{n+1}$  for
which
\begin{equation}  
I_{n+1} \,  \cap \, \Delta(r/q) \, \neq  \, \emptyset   \ \ \mbox{
for some \ } r/q  \in C(n)  \ .
\end{equation}
The sought after collection $ \JJ_{n+1}$ consists precisely  of those intervals that survive.   Formally, for $n \ge 1 $  we let
$$
\JJ_{n+1}:=\{I_{n+1}\in\II_{n+1}\;:\;I_{n+1} \,  \cap \,
\Delta(r/q)=\emptyset\ \mbox{ for any }r/q\in C(n)\}.
$$
For any interval $J_{n-1}   \in \JJ_{n-1}$ and any integer $ R \ge e^{12}$,  we claim that
\begin{eqnarray}
\#\{I_{n+1}\in \II_{n+1}  \!\!\!\! &:&\!\!\!\! J_{n-1}   \cap I_{n+1}\cap
\Delta(r/q)\neq\emptyset\ \mbox{ for some }\ r/q\in C(n)\} \nonumber  \\[2ex]
&\le &
7 \, \log^2\!Rn^2(\log^*\!n)^2  \, . \label{ineq_jn1}
\end{eqnarray}
It then follows from the definition of   $r_{m,n}$  that
$$
\#\{I_{n+1}\in\II_{n+1}\backslash\JJ_{n+1}\;:\; I_{n+1}\subset
J_{n-1}\}\le7\log^2\!R\cdot n^2(\log^*\!n)^2   \stackrel{\eqref{def_rrnm}}=r_{n-1,n}
$$
and therefore the act of removing `bad' intervals from $\II_{n+1}$  is exactly in keeping with the
removal procedure associated with the  construction of a  $({\rm I},\RRR,\vr)$ Cantor set. The goal now is to justify \eqref{ineq_jn1}.

\subsubsection{Counting removed intervals}\label{seq_counting}

\noindent{\em Stage 1.} Let  $ r/q \in C(n)$.  Then there exists a non-negative integer $k$ and an integer $\bar{q}$  such that
\begin{equation}\label{def_qk}
q=D_k\cdot \bar{q}  \qquad\mbox{and   } \quad q\not\in D_{k+1}\ZZ.
\end{equation}
%
%
%
Then,   $$  H(q):=D_k\cdot \bar{q}^{2}   \, . $$
Since all the terms $d_k$ of $\DDD$ are greater than or equal to  two, we have that
\begin{equation}\label{ineq_ek}
D_k \le 2^k \, .
\end{equation}
Next, note that  $q^2\ge H(q)\ge R^{n-1}F(n-1)$. Thus, for any $R > e^9$ it follows  that
\begin{eqnarray}
f(q)  & \ge  &  \mbox{$\frac12$} \, \log^*\! (R^{n-1}F(n-1)) \log^*\!\mbox{$\frac12$}\log
(R^{n-1}F(n-1))   \nonumber \\[2ex]
&\ge& \mbox{$\frac12$}  \, n(\log^*\!n)^2  \, . \label{ineq_f}
\end{eqnarray}
 To see this first observe  that \eqref{ineq_f} for $n=1$ is clearly true.  For $n \ge  2$, by  Stirling formula we have that
$$
R^{n-1}F(n-1) \, \ge \,  R^{n-1}(n-1)!  \, > \, (8n)^n\quad\mbox{ for any  } \ R > e^9   \, .
$$
Therefore the left hand side of \eqref{ineq_f} is bigger than
$$
\mbox{$\frac12$}n\log^*\!(8n)\cdot \log^*\! \left(\mbox{$\frac12$} n \log
(8n)\right)  \; >  \; \mbox{$\frac12$}n(\log^*\!n)^2.
$$

\vspace*{3ex}

\noindent{\em Stage 2.}
We subdivide the collection $C(n)$ of rationals into various `workable' sub-collections.  In the first instance, for any integer $k \ge 0$, let $\cC(n,k) \subset \cC(n)  $ denote the collection of rationals  satisfying  the additional condition \eqref{def_qk}.  Formally,
\begin{equation}\label{def_cnk}
C(n,k):=\left\{r/q \in C(n)\;:\; q\mbox{ satisfies
\eqref{def_qk}}\right\}.
\end{equation}
For any $r/q\in C(n,k)$ we have that
$H(q)=D_k\cdot \bar{q}^{2}$  and  thus in view of~\eqref{zeq2}
and~\eqref{ineq_ek} it follows that
\begin{eqnarray}
0  \, \le \, k  &\le&    [\log_2(R^n F(n)) ]  \, < \, n\log_2R+n\log_2n+n\log_2\log^*\!n \nonumber \\[1ex]
& < & c_2 \,
n\log^*\! n,  \label{ineq_k}
\end{eqnarray}
where $c_2:=(\log R+2)/\log 2$ is an absolute constant independent on
$n$. The upshot is that for fixed $n$ the number of (non-empty) collections  $C(n,k)$
is at most $c_2n\log^*\! n$.

\noindent Next,  for any integer $ l \ge 0$, let $\cC(n,k,l) \subset \cC(n,k)
$ denote the collection of  rationals  satisfying   the
additional condition that
\begin{equation}\label{def_cnkl}
e^lR^{n-1}F(n-1) \, \le  \,
H(q)  \, <  \, e^{l+1}R^{n-1}F(n-1)  \, .
\end{equation}
Formally,
\begin{equation*}\label{def_cnklsv}
C(n,k,l):=\left\{r/q \in C(n,k)\;:\; q\mbox{ satisfies
\eqref{def_cnkl}}\right\}.
\end{equation*}
In view of \eqref{zeq2} we have that
\begin{equation*}
e^l<R n\log^*\!n
\end{equation*}
and thus it follows that
\begin{eqnarray}\label{ineq_l}
0 \ \le \  l & \le &  [\log (Rn\log^*\! n)] \, < \, \log R + 2\log^*\! n   \nonumber \\[1ex]
&<& c_3 \log^*\! n
\end{eqnarray}
where $c_3:=2+\log R$. The upshot is that for fixed $n$ and $k$  the number of (non-empty) collections  $C(n,k,l)$
is at most $c_3 \log^*\! n$. Notice that within any  collection $C(n,k,l)$ we
have extremely tight control on the height.

\vspace*{3ex}

\noindent{\em Stage 3.}
Fix an  interval $J_{n-1}\in \JJ_{n-1}$. Recall, that our goal is establish  \eqref{ineq_jn1}.  This we will do by estimating the quantity
$$
\#\{I_{n+1}\in \II_{n+1}\;:\; J_{n-1}  \cap  I_{n+1}\cap
\Delta(r/q)\neq\emptyset\ \mbox{ for some }\ r/q\in C(n,k,l)\}
$$
and then summing over all possible values of $k$ and $l$. With this in mind,
consider a rational  $r/q\in C(n,k,l)$ and assume that $R> e^9$. Then
\begin{eqnarray}
\#\{I_{n+1}\in \II_{n+1} \!\!\!\!\! & \!\!\!\!\! :  \!\!\!\!\! \; & \!\!\!\!\!  I_{n+1}\cap
\Delta(r/q)\neq\emptyset \}  \ \le \
\displaystyle\frac{|\Delta(r/q)|}{|I_{n+1}|}+2  \nonumber \\[2ex] &=&
\displaystyle\frac{2cR^{n+1}F(n+1)}{c_1f(q)H(q)}+2  \nonumber \\[2ex]  &\stackrel{\eqref{def_cnkl}}{\le} &
\frac{2c R^2 n(n+1) \; [\log^*\!n] \; [\log^*\!(n+1)]}{c_1 f(q) e^l}+2
\nonumber \\[2ex]
&\stackrel{\eqref{ineq_f}}{<} & \displaystyle\frac{8cR^2
(n+1)}{c_1e^l}+2  \, .   \label{iona}
\end{eqnarray}

\noindent Next, consider two rationals $r_1/q_1$, $r_2/q_2 \in
C(n,k,l)$. By definition, there exit integers  $ \bar{q}_1$, $\bar{q}_2 $ so that
$$
q_1=D_k\bar{q}_1\quad\mbox{and}\quad q_2=D_k\bar{q}_2.
$$
Thus  $(q_1,q_2)\ge D_k$ and we have  that

$$
\left|\frac{r_1}{q_1}-\frac{r_2}{q_2}\right| \ge
\frac{1}{D_k\bar{q}_1\bar{q}_2}=(H(q_1)H(q_2))^{-1/2} \; \stackrel{\eqref{def_cnkl}}> \; e^{-l-1}R^{-n+1}F^{-1}(n-1).
$$

\noindent It is easily verified  that $2\, |\Delta(r/q)| $ is less than the right hand side of the above inequality  -- this makes  use of  the fact that $4ec < 1$  which is true courtesy  of \eqref{ineq_c}.
Therefore,
$$
\Delta(r_1/q_1) \cap \Delta(r_2/q_2)  = \emptyset
$$
and it follows that
\begin{eqnarray}
\#\{  r/q\in C(n,k,l)  \!\!\!\! &\!\!\!\! :  \!\!\!\! \;& \!\!\!\! J_{n-1}  \cap
\Delta(r/q)\neq\emptyset \}
\nonumber \\[2ex]
&\le&
2+\frac{|J_{n-1}|}{e^{-l-1}R^{-n+1}F^{-1}(n-1)}\  =  \  2+c_1e^{l+1}.
 \label{ayesha}
\end{eqnarray}

%
%

\noindent The upshot of the  cardinality estimates \eqref{iona} and \eqref{ayesha}  is that
\begin{eqnarray*}
\#\{I_{n+1}\in \II_{n+1} \!\!\!\! &\!\!\!\! :  \!\!\!\! \;& \!\!\!\! J_{n-1}  \cap  I_{n+1}\cap
\Delta(r/q)\neq\emptyset\ \mbox{ for some } r/q\in C(n,k,l)\}
\nonumber \\[2ex]
&\le&
\left(2+c_1e^{l+1}\right)\left(2+\frac{8cR^2(n+1)}{c_1e^l}\right)  \nonumber \\[2ex]
&=  & 4+ 2c_1e^{l+1}+\frac{16cR^2(n+1)}{c_1e^l}+8ecR^2(n+1)  .
\label{ineq_1}
\end{eqnarray*}
By summing  over  $l$ satisfying \eqref{ineq_l} we find that
\begin{eqnarray*}
\#\{I_{n+1}\in \II_{n+1} \!\!\!\! &:& \!\!\!\! J_{n-1}   \cap  I_{n+1}\cap
\Delta(r/q)\neq\emptyset\ \mbox{ for some } r/q\in C(n,k)\}
\\[2ex]
&\le&\sum_{e^l<Rn\log^*\!n}  \!\!\! 2c_1e^{l+1}+ \sum_{l=0}^{c_3\log^*\! n}
\frac{16cR^2(n+1)}{c_1e^l}+c_3\log^*\!n(4+8ec R^2 (n+1)) \\[2ex]
&\le& c_4 n
\log^*\! n
\end{eqnarray*}
where
$$c_4:= 2e^2c_1R+\frac{64c}{c_1}R^2+(\log R+2)(16e \, c \, R^2+4).$$
By summing over  $k$ satisfying  \eqref{ineq_k} we find that
\begin{eqnarray}
\#\{I_{n+1}\in \II_{n+1} \!\!\!\! &:&  \!\!\!\! J_{n-1}   \cap I_{n+1}\cap
\Delta(r/q)\neq\emptyset\ \mbox{ for some } r/q\in C(n)\}  \nonumber \\[2ex]
&\le&  c_2c_4 n^2(\log^*\!n)^2   \, . \label{prop_main}
\end{eqnarray}
In view of \eqref{ineq_c1} and \eqref{ineq_c}, for any  $R >  e^{12}$ the right hand side of
\eqref{prop_main} is bounded by
\begin{equation*}
\left(2+\frac{4(\log R+2)^2}{\log 2}\right)n^2(\log^*\!n)^2  \ <  \ 7\log^2\!R\cdot n^2(\log^*\!n)^2   \, .
\end{equation*}
This establishes   \eqref{ineq_jn1}  as required.

\section{Proof of Proposition \ref{prop2sv}  modulo Theorem~\ref{th_cantor2}}

The proof of Proposition \ref{prop2sv} follows the same structure and ideas as the proof of Proposition~\ref{prop1sv}. In view of this  it is really  only necessary to point out the key differences.

\vspace*{2ex}

Throughout,  $f:\NN\to \RR  :  q \to f(q)=\log^*\!\log q\cdot \log^*\!\log^*\!\log q$ and $\DDD:=\{2^{2^n}\}_{n\in \NN}.
$
Note that by definition
\begin{equation}\label{ineq_dk}
D_k\ge 2^{2^k}  \,  .
\end{equation}

\noindent With $R $, $ c_1$ and  $c$  as in  the proof of Proposition~\ref{prop1sv}, the basic construction of
$$\KK_{\DDD}(f,c) :=   \bigcap_{n=1}^\infty \bigcup_{J\in\JJ_n}J  \  \subset \  \badl_{\DDD}(f,c)$$
 remains pretty much  unchanged  apart from the fact that  the function $F$ is given by
$$
F(n):=\prod_{k=1}^n [\log^*\! k\cdot\log^*\!\log k] \;\  \mbox{ for }
n\ge 1\quad\mbox{ and }  \quad F(0):=1  \   \mbox{ for }n\le 0  .
$$
Also, it becomes apparent from the construction of the collections $\JJ_n$ that
$ \KK_{\DDD}(f,c) $  is a  $({\rm I},\RRR,\vr)$ Cantor set  with   $\RRR=(R_n)$ given by
\begin{equation*}\label{def_rnsv}
R_n:=R \, [\log^*\!(n+1)\log^*\!\log(n+1)]
\end{equation*}
and  $\vr=(r_{m,n})$  given by
\begin{equation*}\label{def_rrnmsv}
r_{m,n} \, :=  \, \left\{\begin{array}{ll} 7 \, \log^2\!R(\log^*\!n)^2(\log^*\!\log n)^2 &\mbox{
\ if }\; m=n-1\\[2ex]
0    &\mbox{ \ otherwise. }
\end{array}\right.
\end{equation*}
Then, it is easily verified that \eqref{cond_th2} is valid for any $R > e^9$ and so  Proposition~\ref{prop2sv}  follows via Theorem~\ref{th_cantor2}.

\vspace*{2ex}

%

Regarding the construction of the collections $\JJ_n$, the  induction procedure is precisely as in \S\ref{cone}. The upshot is that the proof of Proposition~\ref{prop2sv}    reduces to establishing the  following analogue of \eqref{ineq_jn1}. For any interval  $J_{n-1}   \in \JJ_{n-1}$ and any integer $ R > e^{12}$, we have that
\begin{eqnarray}
\#\{I_{n+1}\in \II_{n+1}  \!\!\!\! &:&\!\!\!\! J_{n-1}   \cap I_{n+1}\cap
\Delta(r/q)\neq\emptyset\ \mbox{ for some }\ r/q\in C(n)\} \nonumber  \\[2ex]
&\le &
 7 \, \log^2\!R(\log^*\!n)^2(\log^*\!\log n)^2 \, . \label{ineq_jn1sv}
\end{eqnarray}
This implies that  act of removing `bad' intervals from $\II_{n+1}$  when constructing $\JJ_{n+1} $ from $\JJ_n$ is exactly in keeping with the
removal procedure associated with the  construction of a  $({\rm I},\RRR,\vr)$ Cantor set.
In order to establish \eqref{ineq_jn1sv} we follow the arguments set out in  \S\ref{seq_counting}. For completeness and ease of comparison we briefly describe the analogue of the key estimates.

\vspace*{3ex}

\noindent{\em Stage 1.}  The  analogue of  \eqref{ineq_f}  is the statement that for any  $R >  e^4$
\begin{eqnarray}
f(q)&\ge &\log^*\! \mbox{$\frac12$} \log(R^{n-1}F(n-1))\cdot \log^*\!\log^*\!  \mbox{$\frac12$}
\log(R^{n-1}F(n-1))\nonumber\\[2ex]
&\ge &\log^*\! n\log^*\!\log n.\label{ineq_ff}
\end{eqnarray}
This makes use of the fact that for $n \ge 2 $
$$
R^{n-1}F(n-1)\ge e^{2n} \quad \mbox{ for any  } \ R >e^4.
$$

\vspace*{2ex}

\noindent{\em Stage 2.} In view of \eqref{ineq_dk}, it follows that the  analogue of
\eqref{ineq_k} is that
\begin{equation}\label{ineq_kk}
0 \, \le \, k  \, \le  \,  [\log_2\log_2(R^nF(n))] \, <  \, \tilde{c}_2 \log^*\!n
\end{equation}
where
$$\tilde{c}_2:=\frac{1}{\log 2}\left(2+\log\frac{\log R+2}{\log 2}\right)<c_2  \,  $$
Note that  $ \tilde{c}_2 < c_2 $ is valid since  $R \ge 6$.  Next, in view of \eqref{zeq2} we have that
\begin{equation*}\label{ineq_ell}
e^l \, <  \, R \log^*\!n\log^*\!\log n
\end{equation*}
and thus it follows that
\begin{equation}\label{ineq_ll}
0 \ \le \ l \  \le  \  c_3\log^*\!\log n  \, .
\end{equation}

\vspace*{2ex}

\noindent{\em Stage 3.} Fix an  interval $J_{n-1}\in \JJ_{n-1}$.
Then  \eqref{ayesha}  remains unchanged and the analogue
of~\eqref{iona} is  as follows.  Consider a rational  $r/q\in
C(n,k,l)$ and assume that $R> e^4$. Then
\begin{eqnarray*}
\#\{I_{n+1}\in \II_{n+1} \!\!\!\!\! & \!\!\!\!\! :  \!\!\!\!\! \; & \!\!\!\!\!  I_{n+1}\cap
\Delta(r/q)\neq\emptyset \}  \ \le \
\displaystyle\frac{2cR^{n+1}F(n+1)}{c_1f(q)H(q)}+2  \nonumber \\[2ex]  &\stackrel{\eqref{def_cnkl}}{\le} &\frac{2cR^2 \; [\log^*\!n\log^*\!\log
n]  \; [\log^*\!(n+1)\log^*\!\log(n+1)]}{c_1f(q)e^l}+2
\\[2ex]
&\stackrel{\eqref{ineq_ff}}{<} &
\displaystyle\frac{8cR^2\log(n+1)\log^*\!\log(n+1)}{c_1e^l}+2 \, .
\end{eqnarray*}

\noindent The upshot is that
\begin{eqnarray*}
\#\{I_{n+1}\in \II_{n+1} \!\!\!\! &\!\!\!\! :  \!\!\!\! \;& \!\!\!\! J_{n-1}  \cap  I_{n+1}\cap
\Delta(r/q)\neq\emptyset\ \mbox{ for some } r/q\in C(n,k,l)\}
\nonumber \\[2ex]
&\le&
\left(2+c_1e^{l+1}\right)\left(2+\frac{8cR^2\log^*\!(n+1)\log^*\!\log(n+1)}{c_1e^l}\right).
\end{eqnarray*}
By summing up over  $l$ satisfying \eqref{ineq_ll} we find that
\begin{eqnarray*}
\#\{I_{n+1}\in \II_{n+1}\!\!\!\! &:& \!\!\!\! J_{n-1}   \cap  I_{n+1}\cap
\Delta(r/q)\neq\emptyset\ \mbox{ for some } r/q\in
C(n,k)\}\\[2ex]
&\le & \sum_{e^l<R\log^*\!n\log^*\log n} \!\!\!\!\!\!\!\!\!\!\!\! 2c_1e^{l+1} \ + \
\sum_{l=0}^{c_3\log^*\!\log n} \frac{16cR^2\log^*\!(n+1)\log^*\! \log
(n+1)}{c_1e^l}\\[2ex]
&& ~ \hspace*{6ex} + \ c_3\log^*\!\log n(4+8 e \, c \, R^2 \log^*\!(n+1)\log^*\!\log (n+1))  \\[2ex]
& \le &  c_4
\log^*\!n (\log^*\!\log^*\! n)^2  \, .
\end{eqnarray*}

\noindent By summing up over $k$ satisfying \eqref{ineq_kk} we find that
\begin{eqnarray*}
\#\{I_{n+1}\in \II_{n+1}\!\!\!\!&:& \!\!\!\!  J_{n-1}   \cap I_{n+1}\cap
\Delta(r/q)\neq\emptyset\ \mbox{ for some } r/q\in C(n)\} \nonumber \\[2ex]
& \le  &  c_2c_4 (\log^*\!n)^2\cdot(\log^*\! \log n)^2 \label{prop_main2}  .
\end{eqnarray*}
In view of \eqref{ineq_c1} and \eqref{ineq_c}, for any  $R>e^{12}$   the right hand side of
the above  inequality is bounded by
\begin{eqnarray*}
\left(2+\frac{4(\log R+2)^2}{\log
2}\right)(\log^*\!n)^2\cdot(\log^*\!\log n)^2
&< &7\log^2\!R\cdot (\log^*\!n)^2\cdot(\log^*\!\log n)^2   \, .
\end{eqnarray*}
This establishes  \eqref{ineq_jn1sv} and thereby completes the proof of Proposition~\ref{prop2sv}.

\section{Preliminaries for Theorem~\ref{th_cantor2} \label{Stalin}}

The overall strategy is simple enough. We show that under the hypothesis of the theorem, any given  set $\KK({\rm I},\RRR,\vr)$ contains a `local'  subset $\LKK({\rm I},\RRR,\vs)$  satisfying the desired lower bound inequality for the Hausdorff dimension.   A general and classical method
for obtaining a lower bound for the Hausdorff dimension of an
arbitrary set is the following mass distribution  principle -- see
\cite[pg. 55]{falc}.

\begin{MDP}
Let $ \mu $ be a probability measure supported on a subset $X$ of $
\RR$. Suppose there are  positive constants $a, s $ and $l_0$ such
that
\begin{equation}\label{mdp_eq1}
\mu  ( B ) \le \, a \;   |B|^s \; ,
\end{equation}
for any interval $B$ with length $|B|\le l_0$. Then,  $\dim X \ge s
$.
\end{MDP}

\noindent As we shall soon see, the construction of the local set alluded to above is much simpler than that of  $\KKK({\rm I},\RRR,\vr)$  and enables us to exploit the mass distribution principle.

\subsection{Local Cantor sets} A  $({\rm I},\RRR,\vr)$ Cantor set $\KKK({\rm I},\RRR,\vr)$ is said to be  local if $r_{m,n}=0$ whenever  $m\neq n$. Furthermore, we write  $\LKK({\rm I},\RRR,\vs)$  for $\KKK({\rm I},\RRR,\vr)$  where
$$
\vs:=(s_n)_{n\in \ZZ_{\ge 0}}  \quad   {\rm and } \quad
s_n:=r_{n,n}.
$$
The set $\LKK({\rm I},\RRR,\vs)$ will be referred to as a {\em $({\rm I},\RRR,\vs)$ local Cantor set.}

\vspace*{2ex}

In a nutshell, the removing procedure associated with the construction of a local Cantor set has no `memory' -- it depends only on the level under consideration.    More formally,   given the collection $\JJ_{n}$ of level $n$ survivors,  the construction of $\JJ_{n+1} $  is completely independent of the previous level $k$ ($<n$) survivors.  Indeed the  construction is totally   local  within  each  interval $J_{n} \in \JJ_{n}$.  It is this fact that is utilized when attempting to establish the following dimension  result  for the  associated local Cantor set.  Note that in view of Theorem~\ref{th_cantor1},  any local set  $\LKK({\rm I},\RRR,\vs)$ is non-empty if
$
R_n-s_n>0 $ for all $ n\in \ZZ_{\ge 0}
$.

\begin{lemma}\label{lem_local}
Given   $\LKK({\rm I},\RRR,\vs)$, suppose that
$$
t_n:=R_n-s_n   > 0  \quad \forall \  n\in\ZZ_{\ge 0}  \, .
$$
Furthermore, suppose there are positive constants  $s$ and $n_0 $
such that for all  $n>n_0$
\begin{equation}\label{ineq_lemloc}
R_n^s\le t_n   \, .
\end{equation}
Then
$$
\dim \LKK({\rm I},\RRR,\vs) \ge s.
$$
\end{lemma}

\noindent{\em Proof. \, }  We start by constructing  a  probability measure
$\mu$ supported on $$\LKK({\rm I},\RRR,\vs) :=   \bigcap_{n=1}^\infty \bigcup_{J\in
\JJ_n} J.$$ in the standard manner. For any $J_n \in
\JJ_n$, we attach a weight $\mu(J_n)$ defined recursively as
follows.

\vspace*{2ex}

\noindent For  $n=0$, $$ \mu(J_0)\ := \frac{1}{\#\JJ_0}=1  \ $$ and
for $n\ge 1$,
\begin{equation}\label{beq4}
\mu(J_n)  \, :=  \, \frac{\mu(J_{n-1})}{\# \{J\in \JJ_n\;:\;
J\subset J_{n-1}\}} \
\end{equation}

\vspace*{2ex}

\noindent where $J_{n-1} \in \JJ_{n-1}$   is the unique interval
such that $J_n\subset J_{n-1}$.
 This procedure thus defines inductively a mass on any interval
appearing in the construction of $\LKK({\rm I},\RRR,\vs) $. In fact a lot more is true
---   $\mu$  can be further extended to all Borel subsets $F$ of
$\RR $  to determine $\mu(F) $  so that $\mu$ constructed as above
actually defines  a measure supported on $\LKK({\rm I},\RRR,\vs) $. We now state this
formally.

\begin{itemize}
\item[]  {\em Fact.} The probability measure $\mu$ constructed above is
supported on $\LKK({\rm I},\RRR,\vs) $ and for any Borel set $F$
\[
\mu(F):= \mu(F \cap\LKK({\rm I},\RRR,\vs) )  \; = \; \inf\;\sum_{J\in  \cJ }\mu(J)  \
.
\]
The infimum is over all coverings $\cJ$ of $F \cap \LKK({\rm I},\RRR,\vs) $ by
intervals  $J\in \{\JJ_n: n \in \ZZ_{\ge 0} \}$.
\end{itemize}

\noindent For further details see \cite[Prop. 1.7]{falc}.   It remains to show that $\mu$ satisfies \eqref{mdp_eq1}. Firstly,  notice that  for any interval $J_n
\in \JJ_n$ we have that
\begin{equation}\label{svdelta}
\mu(J_n) \stackrel{\eqref{beq4}}\le t_{n-1}^{-1}\ \mu(J_{n-1}) \le \prod_{i=0}^{n-1}
t_i^{-1}
\end{equation}
Next, let $\delta_n$ denote the length of a generic interval $J_n
\in \JJ_n$. In view of the splitting procedure associated with the construction of $\LKK({\rm I},\RRR,\vs)$,  we find that
\begin{equation}\label{eq_delta}
\delta_n=|I|\cdot \prod_{i=0}^{n-1}R_i^{-1} \, .
\end{equation}
Consider an arbitrary interval $B\subset [0,1] $ with length $|B| <
\delta_{n_0}$. Then there exists an integer $ n\ge n_0 $ such that
\begin{equation} \label{fine}
\delta_{n+1}   \, \le  \, |B|  \,  <  \, \delta_n  \; .
\end{equation}
It follows that
\begin{eqnarray*}
\mu(B)  &\le& \sum_{\substack{J \in \JJ_{n+1}:\\
\scriptstyle J \cap B \neq\emptyset }}
\!\! \mu(J)    \ \ \stackrel{\eqref{svdelta}}\le \ \   \left\lceil
\frac{|B|}{\delta_{n+1}}\right\rceil
\prod_{i=0}^n t_i^{-1} \\[2ex]
& \stackrel{\eqref{eq_delta}}\le  &  2 \,
\frac{|B|}{|I|}\prod_{i=0}^n \frac{R_i}{t_i}   \ \ = \ \  2 \,
\frac{|B|}{|I|}^{1-s}  \  \ \prod_{i=0}^n
\frac{R_i}{t_i}  \, \cdot \ |B|^s \\[2ex]
& \stackrel{\eqref{fine}}<  & 2 \,
\frac{\delta_n}{|I|}^{\!\!\!\!1-s}  \ \   \prod_{i=0}^n
\frac{R_i}{t_i}  \, \cdot \ |B|^s  \\[2ex]
& \stackrel{\eqref{eq_delta}}<  &   2  |I|^{-s}\prod_{i=0}^n
\frac{R_i^{s}}{t_i} \ \cdot  \
|B|^{s}   \\[2ex]
&  \stackrel{\eqref{ineq_lemloc}}{\le}  & 2 \,
|I|^{-s}\prod_{i=0}^{n_0}\frac{R_i^{s}}{t_i}  \ \cdot \ |B|^{s}   \,
.
\end{eqnarray*}
In other words,   \eqref{mdp_eq1} is valid with
$$a:=2 \,
|I|^{-s}\prod_{i=0}^{n_0}\frac{R_i^{s}}{t_i}  \, $$ and on applying
the mass distribution principle we obtain the desired statement.
\newline\hspace*{\fill}$\boxtimes$

\vspace*{2ex}


 In view of  Lemma~\ref{lem_local},  the proof of Theorem~\ref{th_cantor2} reduces to establishing the following key statement.

\begin{proposition}\label{lem_subcantor}
Let  $\KKK({\rm I},\RRR,\vr)$   be as in Theorem~\ref{th_cantor2}.  Then
there exists a local Cantor type set  $$\LKK({\rm I},\RRR,\vs)  \subset
\KKK({\rm I},\RRR,\vr) $$ where
$$
\vs:=(s_n)_{n\in \ZZ_{\ge 0}}  \quad   { with  } \quad  s_n:=
\mbox{$\frac12$} \,  R_n \, .
$$
\end{proposition}

\vspace*{2ex}

\noindent Indeed, by Proposition~\ref{lem_subcantor} we have that
$$
\dim \KKK({\rm I},\RRR,\vr)\ge \dim  \LKK({\rm I},\RRR,\vs) .
$$
Now fix some positive $s<\liminf\limits_{n\to\infty}(1-\log_{R_n}\!2)$. Then,  there
exists an integer  $n_0$  such that
$$
s \, <  \, 1-\log_{R_n}\!2  \quad {\rm \ for \  all }  \quad  n>n_0 \, .
$$
Also  note that
$$
t_n=R_n-s_n =\frac{R_n}{2}
$$
and
$$
R_n^s<\frac{R_n}{2}= t_n   \quad {\rm \ for \  all }  \quad  n>n_0
\, .
$$
Therefore, Lemma~\ref{lem_local} implies that
$$\dim \LKK({\rm I},\RRR,\vs)  \ge s   \, . $$ This inequality is true for any
$s<\liminf\limits_{n\to\infty}(1-\log_{R_n}\!2)$ and  hence  completes
the proof of Theorem~\ref{th_cantor2} modulo Proposition~\ref{lem_subcantor}.

\vspace*{2ex}

Before moving on to the proof of the proposition, it is useful to first
investigate  the distribution of intervals within each collection $\JJ_n$ associated with $\KKK({\rm I},\RRR,\vr)$.


\subsection{The distribution of intervals within $\JJ_n$}

In this section,  the set
$$
 \KKK ({\rm I},\RRR,\vr) :=  \bigcap_{n=1}^\infty \bigcup_{J\in
\JJ_n} J
$$
 and the sequence $\vs$ are  as in Proposition~\ref{lem_subcantor}.

Let $\TTT_0:= \{ {\rm I}\}$. For $n \ge 1$,
let $\TTT_n$ denote a generic  collection of intervals obtained  from
$\TTT_{n-1}$ via the splitting and removing procedures
associated with a  $({\rm I},\RRR,\RRR-\vs)$ local Cantor set.  Here $\RRR-\vs$ is the sequence
$(R_n-s_n)$. Then, clearly
$$
\# \TTT_{n+1}   \ \ge  \    \# \TTT_{n}  \times  s_n  \,   \quad \forall \ n \in \ZZ_{\ge0} \, .
$$
Loosely speaking, the following result shows that the intervals
$J_n$ from $\JJ_n$ are ubiquitous within the interval ${\rm I}$.

\begin{lemma}\label{sl_lem1} For $R$ sufficiently large,
\begin{equation}\label{sl_lem_eq}
\TTT_n\cap \JJ_n \neq \emptyset   \qquad \forall \ \ n \in \ZZ_{\ge0} \, .
\end{equation}
\end{lemma}

\noindent{\em Proof. \, } For an integer $n \ge 0$, let $h(n)$
denote the cardinality of the set $\TTT_n\cap \JJ_n$. Trivially,
$h(0)=1$ and lemma would follow on showing that
\begin{equation}\label{goodcount}
h(n+1)\ge \frac{R_n}{4} \,  h(n)  \, .
\end{equation}



\noindent  for all $ n \in \ZZ_{\ge0}$.  This we now do via induction.  Consider the set $\TTT_n\cap \JJ_n$. By the construction of $\TTT_{n+1}$ and the splitting procedure associated with $\KKK({\rm I},\RRR,\vr)$,  each
of the $h(n)$ intervals in $\TTT_n\cap \JJ_n$ gives rise to at least
$s_n$ intervals $I_{n+1}$ in $\TTT_{n+1}\cap \II_{n+1}$. By the
the removing procedure associated with $\KKK({\rm I},\RRR,\vr)$,
for each interval  $J_k\in\TTT_k\cap\JJ_k$ ($ 0 \le k \le n$)  we remove  at
most $r_{k,n}$ intervals $I_{n+1}\in \TTT_{n+1}\cap \II_{n+1}$ that lie within
$J_k$. The upshot of this is that
\begin{equation}\label{ineq_hn}
h(n+1)\ge s_n h(n)-\sum_{k=0}^{n} r_{k,n}\, h(k).
\end{equation}
For $n=0$ this inequality is  transformed to
$$
h(1)\ge \frac{R_0}{2}h(0)-r_{0,0}h(0)\ \stackrel{\eqref{cond_th2}}\ge\  \frac{R_0}{4}h(0)
$$
as required.  Now assume
that \eqref{goodcount} is valid for all $ 1 \le m \le n$. In
particular,  it means that
$$
h(m)\le \frac{4}{R_m} \, h(m+1)\le\ldots\le
\prod_{k=1}^{n-m}\frac{4}{R_{n-k}}\, h(n)
$$
which together with   \eqref{ineq_hn}  implies that
\begin{eqnarray*}
h(n+1)  & \ge  &  \frac{R_n}{2}h(n)-\left(\sum_{k=0}^n
\left(r_{n-k,n}\prod_{i=1}^k\frac{4}{R_{n-i}}\right)\right)h(n)  \\[2ex]
& \stackrel{\eqref{cond_th2}}\ge &   \frac{R_n}{4}  \, h(n)  \, .
\end{eqnarray*}
This completes the induction step and  thus establishes  the desired inequality
\eqref{goodcount}  for all $n \in \ZZ_{\ge 0}$.
\hspace*{\fill}$\boxtimes$

\section{Proof of Proposition \ref{lem_subcantor} \label{Lenin}}

 By definition, the set $\KKK({\rm I},\RRR,\vr)$ is the intersection of closed intervals $J_n$ lying within nested collections $\JJ_n$.
 For
each integer $n \ge 0$, the aim is to construct a nested collection $\LLL_n
\subseteq \JJ_n$ that  complies with the  construction of  a   $({\rm I},\RRR,\vs)$ local Cantor set. Then, it  would follow that
$$
\bigcap_{n=0}^\infty \bigcup_{J\in
\LLL_n} J
$$
is precisely the desired set $\LKK({\rm I},\RRR,\vs)$.

\subsection{Construction the collection  $\LLL_n$}

For any integer $n \ge 0$, the goal of this section is to construct the
desired nested collection $\LLL_n \subseteq \JJ_n$ alluded to above.  This will involve constructing auxiliary collections
$\LLL_{m,n}$ and $\R_{m,n}$ for integers $m,n$ satisfying $0\le m\le
n$. For a fixed $n$,  let
$$
\JJ_0  \, ,  \ \JJ_1 \, ,  \  \ldots, \ \JJ_n   \;
$$
be the collections arising from  the construction of $\KKK({\rm
I},\RRR,\vr)$. We will require $\LLL_{m,n}$ to satisfy the following
conditions.
\begin{itemize}
\item[\bf C1.] For any $0\le m\le n$, $\LLL_{m,n}\subseteq \JJ_m$.
\item[\bf C2.] For any $0\le m< n$, the  collections $\LLL_{m,n}$ are nested; that is
    $$\bigcup_{J\in \LLL_{m+1,n}}J  \qquad \subset \quad \bigcup_{J\in \LLL_{m,n}}J.$$
\item[\bf C3.] For any $0\le m<n$ and $ J_m\in \LLL_{m,n}$, there are at
least  $R_m-s_m$  intervals $J_{m+1}\in \LLL_{m+1,n}$ contained
within $J_m$; that is
\begin{equation*}
\# \{J_{m+1} \in \LLL_{m+1,n} \;:\; J_{m+1} \subset J_{m}\}   \ \ge
\ R_m-s_m \ .
\end{equation*}
\end{itemize}

\noindent  In addition, define $\R_{0,0} := \emptyset $ and for $n
\ge 1$
\begin{equation}\label{def_Rnn}
\R_{n,n}:=\left\{I_{n}\in \II_n\backslash \JJ_n\;:\; I_n\subset
J_{n-1}\mbox{ for some }J_{n-1}\in \LLL_{n-1,n-1}\right\}   \ .
\end{equation}
Furthermore, for $ 0\le m< n $ define
\begin{equation}\label{def_rnm}
\R_{m,n}:=\R_{m,n-1}\cup \{J_m\in \LLL_{m,n-1}\;:\; \#\{J_{m+1}\in
\R_{m+1,n}\;:\; J_{m+1}\subset J_m\}\ge s_m\ \} \ .
\end{equation}
Loosely speaking and  with reference to condition (C3),  the
collections $\R_{m,n}$ are the `dumping ground'  for those intervals
$J_m\in\LLL_{m,n-1}$ which do not contain enough sub-intervals
$J_{m+1}$. Note that for $n$ fixed, the collections $\R_{m,n}$ are
defined  in descending order with respect to~$m$. In other words, we
start with $\R_{n,n}$  and finish with $\R_{0,n}$.

\vspace*{1ex}

The construction is as follows.

\vspace*{1ex}

\noindent{\em Stage 1.} Let $\LLL_{0,0}:=\JJ_0  $  and   $
\R_{0,0}:=\emptyset$.

\vspace*{1ex}

\noindent {\em Stage 2.} Let $  0 \le t \le n $. Suppose  we have
constructed the desired collections
$$\LLL_{0,t}\subseteq \JJ_0, \   \LLL_{1,t}\subseteq \JJ_1,\ldots, \LLL_{t,t}\subseteq
\JJ_t$$ and
$$\R_{0,t},\ldots, \R_{t,t} \, .  $$

We  now construct the corresponding  collections for $t=n+1$.

\noindent{\em Stage 3.} Define
$$
\LLL'_{n+1,n+1}:=\{J_{n+1}\in \JJ_{n+1}\;:\; J_{n+1}\subset
J_n\mbox{ for some }J_n\in \LLL_{n,n}\}
$$
and let $\R_{n+1,n+1}$ be given by \eqref{def_Rnn} with $n+1$
instead of $n$. Thus the collection  $\LLL'_{n+1,n+1}$ consists of
`good' intervals from $\JJ_{n+1}$ that are contained within  some
interval from $\LLL_{n,n}$.  Our immediate task is to construct the
corresponding collections $ \LLL'_{u,n+1}$ for each $ 0 \le u \le
n$.  These will be constructed  together with the `complementary'
collections $\R_{u,n+1}$ in descending order with respect to $u$.

\vspace*{1ex}

\noindent{\em Stage 4.} With reference to Stage 3, suppose we have
constructed the collections $\LLL'_{u+1,n+1}$ and $\R_{u+1,n+1}$ for
some $0\le u\le n$. We now construct  $\LLL'_{u,n+1}$ and
$\R_{u,n+1}$. Consider the collections $\LLL_{u,n}$ and $\R_{u,n}$.
Observe that some of the intervals $J_u$ from $\LLL_{u,n}$ may
contain less than $R_u-s_u$ sub-intervals from $\LLL'_{u+1,n+1}$ (or
in other words, at least $s_u$ intervals from $\R_{u+1,n+1}$). Such
intervals $J_u$ fail the counting condition (C3)  for $\LLL_{u,n+1}$
and informally speaking are moved out of $\LLL_{u,n}$ and into
$\R_{u,n}$.  The resulting sub-collections are $\LLL'_{u,n+1}$ and
$\R_{u,n+1}$ respectively. Formally,
$$
\LLL'_{u,n+1}:=\{J_u\in \LLL_{u,n}\;:\; \#\{J_{u+1}\in
\R_{u+1,n+1}\;:\; J_{u+1}\subset J_u\}< s_u\ \} \
$$
and $\R_{u,n+1}$ is given by \eqref{def_rnm} with $m$ replaced  by $u$ and $
n$ replaced by $n+1$.

\vspace*{1ex}

\noindent{\em Stage 5.} By construction the collections
$\LLL'_{u,n+1}$ satisfy conditions (C1) and (C3).  However,  for
some $J_{u+1}\in \LLL'_{u+1,n+1}$ it  may be the case that $J_{u+1}$
is not contained in any interval $J_u\in \LLL'_{u,n+1}$  and thus
the collections $\LLL'_{u,n+1}$ are not necessarily nested. The
point is that during  Stage 4 above the interval $J_u\in \JJ_u$
containing $J_{u+1}$ may be  `moved' into $\R_{u,n+1}$. In order to
guarantee the nested condition (C2)  such intervals $J_{u+1}$ are
removed from $\LLL'_{u+1,n+1}$. The resulting sub-collection is the
required auxiliary collection $\LLL_{u+1,n+1}$. Note that
$\LLL_{u+1,n+1}$ is constructed via $\LLL'_{u+1,n+1}$ in ascending
order with respect to $u$.  Formally,
$$
\LLL_{0,n+1}:=\LLL'_{0,n+1}
$$
and for $1\le u\le n+1$
$$
\LLL_{u,n+1}:=\{J_u\in \LLL'_{u,n+1}\;:\; J_u\subset J_{u-1} \mbox{
for some } J_{u-1}\in \LLL_{u-1,n+1}\}   \, .
$$
With reference to Stage 2, this completes the  induction step and
thereby   the construction of the auxiliary collections.

For any integer $n \ge 0$, it remains to  construct the sought
after collection   $\LLL_n$ via the auxiliary collections
$\LLL_{m,n}$. Observe that since
$$
\LLL_{m,m}\supset \LLL_{m,m+1}\supset \LLL_{m,m+2}\supset \ldots
$$
and the cardinality of each collection $\LLL_{m,n}$ with    $m\le n$
is finite, there exists some  integer $N(m)$ such that
$$
\LLL_{m,n}   \ = \ \LLL_{m,n'}  \qquad \forall \quad n,n'  \ge N(m)
\ .
$$

\noindent  Now simply define $$ \LLL_n:=\LLL_{n,N(n)} \ . $$

\noindent Unfortunately, there remains one slight issue.  The
collection   $\LLL_n$  defined in this manner could be empty.

\vspace*{2ex}

The goal now  is to show  that  $\LLL_{m,n}\neq \emptyset$ for any
$m\le n$. This clearly  implies  that $\LLL_n\neq\emptyset$ and
thereby completes the construction.

\subsection{The collection $\LLL_{m,n}$ is non-empty}

\begin{lemma}\label{lem_mnonempty}
For any $m,n\in \NN, m\le n$, the set $\LLL_{m,n}$ is nonempty.
\end{lemma}

\noindent{\em Proof.} Suppose the contrary:
$\LLL_{m,n}=\emptyset$ for some integers  satisfying $ 0 \le m\le
n$. In view of  the construction of $\LLL_{m,n}$ every interval in
$\LLL_{m-1,n}$ contains at least $R_{m-1}-s_{m-1}>0$ sub-intervals
from $\LLL_{m,n}$. Therefore each of the collections $\LLL_{m-1,n},
\LLL_{m-2,n}, \ldots, \LLL_{0,n}$ are also empty and it follows that
$\R_{0,n}=\JJ_0$.

Now consider the set $\R_{m,n}$. By the construction we have the
chain of nested sets
$$
\R_{m,n}\supseteq \R_{m,n-1}\supseteq\cdots\supseteq\R_{m,m}
$$
and in view of \eqref{def_Rnn} the elements of $\R_{m,m}$ are
intervals from $\II_m\backslash \JJ_m$. Consider any interval
$J_m\in \R_{m,n}\backslash \R_{m,m}$. Take $m < n_0 \le n $ such
that $J_m\in \R_{m,n_0}$ but $J_m\not\in \R_{m,n_0-1}$. Then $J_m$
was added to $\R_{m,n_0}$ on stage~4 of the construction. Hence
$I_m$ should have at least $s_m$ sub-intervals  from $\R_{m+1,n_0}$
and therefore from $\R_{m+1,n}$. The upshot of this is the
following: for any interval $I_m$ from $\R_{m,n}$ either $I_m\in
\II_m\backslash \JJ_m$ or $I_m$ contains at least $s_m$
sub-intervals $I_{m+1}\in \R_{m+1,n}$.

Next we exploit Lemma \ref{sl_lem1}. Choose an interval $J_0$ from $
\R_{0,n} = \JJ_0$   and define $\TTT_0:=\{J_0\}$. For $ 0\le m<n$,
we define inductively nested collections
$$
\TTT_{m+1}:=\{I_{m+1}\in \TTT(I_m)\;:\; I_m\in \TTT_m\} 
$$
with $\TTT(I_m)$ given  by one of the following three scenarios.
\begin{itemize}
\item $I_m\in \R_{m,n}$ and $I_m$ contains at least $s_m$ sub-intervals
$I_{m+1}$ from $\R_{m+1,n}$. Let $\TTT(I_m)$ be the collection
consisting of these sub-intervals. Note that when $m=n-1$ we have
$\TTT(I_m)\subset \R_{n,n}\subset \II_n\backslash \JJ_n$. Therefore
$\TTT(I_{n-1})\cap \JJ_{n}=\emptyset$.

\item $I_m\in \R_{m,n}$ and $I_m$ contains strictly less than $s_m$ sub-intervals $I_{m+1}$ from
$\R_{m+1,n}$. Then the interval $I_m \in \II_m\backslash \JJ_m$ and
we  subdivide $I_m$ into $R_m$ closed intervals $I_{m+1}$ of equal
length.  Let $\TTT(I_m)$ be the collection consisting of all of
these sub-intervals. Note that $\TTT(I_m)\cap \JJ_{m+1}=~\emptyset$.

\item $I_m\not\in \R_{m,n}$. Then the interval $I_m$ does not
intersect any interval from $\JJ_m$ and we  subdivide $I_m$ into
$R_m$ closed intervals $I_{m+1}$ of equal length.  Let $\TTT(I_m)$
be any collection consisting of all such sub-intervals. Note that
$\TTT(I_m)\cap \JJ_{m+1}=~\emptyset$.
\end{itemize}

\noindent The upshot is that
$$
\# \TTT_{m+1}   \ \ge \    \# \TTT_m  \times s_m\qquad \forall \; 0<
m\le n
$$
and that
$$
\TTT_{n}\cap \JJ_{n} \, =  \, \emptyset  \; .
$$
However, in view of Lemma~\ref{sl_lem1} the latter is impossible and
therefore the starting premise that $\LLL_{m,n}=\emptyset$ is false.
This completes the proof of Lemma~\ref{lem_mnonempty} and therefore
Proposition~\ref{lem_subcantor}.
\newline\hspace*{\fill}$\boxtimes$

\vspace*{4ex}

\section{Intersecting  Cantor sets  \label{FC}}

With reference to \S\ref{gcf}, fix the interval ${\rm I}$ and the sequence  $\RRR:= (R_n)$. Let $ k  \in \NN$ and consider the two parameter sequences
$$\vr_i:=(r^{\!(i)}_{m,n})   \qquad  1\le i \le k \, .
$$
The following result shows that
the intersection of any finite number of $({\rm I},\RRR,\vr_i) $ Cantor sets is yet another   $({\rm I},\RRR,\vr) $ Cantor set.

\begin{theorem}\label{th_icantor}
For each integer $1\le i \le k $, suppose we are given a set $\KKK
({\rm I},\RRR,\vr_i) $. Then
$$
\bigcap_{i=1}^{k} \KKK ({\rm I},\RRR,\vr_i)   
$$ %
is a $({\rm I},\RRR,\vr) $ Cantor set
where
$$
\vr:=(r_{m,n})  \quad\mbox{with } \quad r_{m,n} :=\sum_{i=1}^k
r^{(i)}_{m,n}  \, .
$$
\end{theorem}

\vspace*{2ex}

\noindent{\em Proof.\, }
Loosely speaking we need to show that there exists a  $({\rm I},\RRR,\vr) $ Cantor set that  simultaneously incorporates the splitting and removing procedures associated with the sets
$$\KKK({\rm I},\RRR,\vr_i) :=  \bigcap_{n=1}^\infty \bigcup_{J\in\JJ^{(i)}_n} J  \qquad \quad   (1\le i \le k) \, . $$
For each $n\in\ZZ_{\ge 0}$, consider the collection
$$
\JJ_n:=\bigcap_{i=1}^k \JJ_n^{(i)}  \, .
$$
We claim  that $\JJ_n$ complies with the construction of a $({\rm I},\RRR,\vr) $ Cantor set.
If true,  then we are done since
$$
\KKK({\rm I},\RRR,\vr):=\bigcap_{n=1}^\infty \bigcup_{J\in\JJ_n} J  = \bigcap_{n=1}^\infty \bigcap_{i=1}^k \bigcup_{J\in\JJ^{(i)}_n} J= \bigcap_{i=1}^k \bigcap_{n=1}^\infty \bigcup_{J\in\JJ^{(i)}_n}J:=\bigcap_{i=1}^{k} \KKK ({\rm I},\RRR,\vr_i)  \, .
$$
Firstly note that the claim is true for $n=0$ since $\JJ_0:=\{{\rm I}\}$.
Now assume that the claim is true for some fixed $n\in\ZZ_{\ge 0}$. Consider an arbitrary interval
$J_n\in \JJ_n$. By definition, $J_n\in \JJ_n^{(i)}$ for each $i$.
By construction,   every interval in $\JJ_n^{(i)}$ gives rise to $R_n$ intervals $I_{n+1}\in \II_{n+1}^{(i)}$. Thus,  for each $J_n\in \JJ_n$ the
collection
$$\II_{n+1}:=\bigcap_{i=1}^k\II_{n+1}^{(i)}
$$
contains exactly $R_n$ intervals $I_{n+1}$ that lie within $J_n$. This coincides precisely with the splitting procedure associated with a $({\rm I},\RRR,\vr) $ Cantor set.   We now turn our attention to the  removing procedure. By construction, for each interval $J_n\in \JJ_{n}^{(i)}$ we remove at most $r_{n,n}^{(i)}$ intervals $I_{n+1}\in\II^{(i)}_{n+1}$ that lie within $J_n$.
Thus for any
$J_n\in\JJ_n$ there are at most
$$
r_{n,n} := \sum_{i=1}^k r_{n,n}^{(i)}
$$
intervals $I_{n+1}\subset J_n$ that are removed from $\II_{n+1}$.
In general, for each $0\le m \le n$ and each interval
$J_{m}\in \JJ_{m}$ there are at most
$$
r_{m,n}:=\sum_{i=1}^k r_{m,n}^{(i)}
$$
additional intervals $I_{n+1}\subset J_m$  that are removed from $ \II_{n+1}$. This coincides precisely with the removing procedure associated with a $({\rm I},\RRR,\vr) $ Cantor set.
The upshot is that $\JJ_{n+1}$ complies with the construction of a $({\rm I},\RRR,\vr) $ Cantor set.  This completes the induction step  and thereby completes the proof of Theorem~\ref{th_icantor}.
\newline\hspace*{\fill}$\boxtimes$

\vspace*{4ex}

\noindent{$\bullet $ \em An application.}
We now describe  a simple application of Theorem~\ref{th_icantor} which enables us to deduce a  non-trivial strengthening of Theorem~\ref{th_main}.
In the course of establishing   Proposition~\ref{prop1sv} we show that the
set $\mad_\DDD(f)$ contains the Cantor-type set $\KKK({\rm I},\RRR,\vr)$
where $\RRR=(R_n)$ and $\vr=(r_{m,n})$ are given by~\eqref{def_rn} and
\eqref{def_rrnm} respectively; namely, for any fixed integer $R > e^{12}$
$$
R_n:=R  \, (n+1)  \, [\log^*\!(n+1)]  \quad {\rm and }  \quad
r_{m,n} \, :=  \, 7\, \log^2\!R\cdot
n^2(\log^*\!n)^2
$$
if $ m=n-1$ and zero otherwise.  Note that although  these quantities are dependent on the actual value of  $R$  the statement that $\KKK({\rm I},\RRR,\vr) \subset \mad_\DDD(f)$  is  not.

 Now for each  $1\le i \le k $, let $\DDD_i$ be a sequence of integers greater than or equal to two and let $f$ be as in  Proposition~\ref{prop1sv}.  Then, with $\RRR$ and $\vr$ as above,   Theorem~\ref{th_icantor} implies that
$$
\bigcap_{i=1}^k \mad_{\DDD_i}(f)\supset
\KKK({\rm I},\RRR,k\vr) \qquad\mbox{ where }  \quad  k\, (r_{m,n}):=(kr_{m,n}).
$$
It is easily verified that for $ R > k \, e^9$
\begin{eqnarray*}
{\rm l.h.s. \ of \ } \eqref{cond_th2}  & = &  k \cdot r_{n-1,n}\cdot\frac{4}{R_{n-1}}  \, \le  \,  k \cdot  7 \cdot  2^3  \cdot \frac{\log^2\!R \cdot n \, \log^*\!n }{R}  \\[2ex]
&\le  &      \frac{ k \cdot  7 \cdot  2^6  \log^2\!R}{R^2} \cdot  \frac{R \, (n+1) \, [\log^*\!(n+1)] }{4} \\[2ex]
& \le &  \frac{R_n  }{4}  \ =  \  {\rm r.h.s. \ of \ }  \eqref{cond_th2}  \, .
\end{eqnarray*}
Hence, for any fixed $ R > k \, e^{12}  $,    Theorem~\ref{th_cantor2}  implies that
$$
\dim\left(\textstyle{\bigcap_{i=1}^k} \mad_{\DDD_i}(f)\right)\ge \liminf_{n\to
\infty}(1-\log_{R_n}\!2)=1.
$$
The complementary upper bound  inequality  for the dimension is trivial. Thus we have established  the
following strengthening of Theorem~\ref{th_main}.

\begin{theorem}\label{thm1sv}  For each  $1\le i \le k $, let $\DDD_i$ be a sequence of integers greater than or equal to $2$  and let $f$ be as in   Proposition~\ref{prop1sv}.
Then
\begin{equation*}
\dim\left(\bigcap_{i=1}^k \mad_{\DDD_i}(f)\right) =1  \, .
\end{equation*}
\end{theorem}

\vspace*{3ex}

\noindent{$\bullet$ \em What about other intersections?} There are two natural problems that arise in relation to Theorem~\ref{th_icantor}.  Firstly, to generalise the statement so as to incorporate any finite number of sequences $\RRR_i:=(R_n^{(i)})$.    Secondly,  to establish the  analogue of
Theorem~\ref{th_icantor} for countable intersections. This is more challenging  than the first and  in all likelihood will involve imposing  extra conditions on the sequences $\RRR$ and $\vr$.
A direct consequence of the `correct' countable version of Theorem~\ref{th_icantor} would be the statement that
$$
\dim\left(\textstyle{\bigcap_{i=1}^{\infty}} \mad_{\DDD_i}(f)\right) =  1 \, .
$$
Note that establishing  the countable analogue  of Theorem~\ref{thm1sv} remains an open problem.

\vspace*{3ex}

\noindent{$\bullet$ \em A more general Cantor framework.}  The Cantor framework of  \S\ref{gcf} and indeed of this section is  one-dimensional. Naturally it would be interesting to develop the  analogous   $n$--dimensional Cantor framework in which intervals are replaced by balls. Establishing the higher dimensional generalisation   of Theorem~\ref{th_cantor2}  and indeed Theorem~\ref{th_icantor} will almost certainly make use of standard covering arguments  from geometric measure theory; for example, the `$5r$' and Besicovitch covering lemmas.    Beyond higher dimensions,   it would be highly desirable to develop an analogue of the framework of  \S\ref{gcf} within the context of `reasonable' metric spaces  --  such as  a (locally) compact metric space equipped with an Ahlfors regular measure. A generalisation of this type would  enhance the scope of  potential applications.

\vspace*{8ex}

\noindent{\bf Acknowledgements.} SV would like to thank all those at Talbot Primary School who have made the dynamic  duo most welcome! In particular a special thanks to Julia Alvin,  David Young and Paulin Jacobson. Also many thanks to Victor Beresnevich and  Maurice Dodson for various suggestions that have improved the clarity  of the paper and to Andrew Pollington for many discussions over many many  years regarding Littlewood.


\def\cprime{$'$}

\vspace{5mm}

\noindent Dzmitry A. Badziahin: Department of Mathematics,
University of York,

\vspace{0mm}

\noindent\phantom{Dzmitry A. Badziahin: }Heslington, York, YO10 5DD,
England.


\noindent\phantom{Dzmitry A. Badziahin: }e-mail: db528@york.ac.uk

\vspace{5mm}

\noindent Sanju L. Velani: Department of Mathematics, University of
York,

\vspace{0mm}

\noindent\phantom{Sanju L. Velani: }Heslington, York, YO10 5DD,
England.


\noindent\phantom{Sanju L. Velani: }e-mail: slv3@york.ac.uk

\end{document}